\documentclass[reqno,11pt]{amsart}
\usepackage{amssymb}

\usepackage{amsmath}
\usepackage{times}

\newtheorem{theorem}{Theorem}[section]

\newtheorem{proposition}[theorem]{Proposition}

\begin{document}
\title{Bounds for mixing time of quantum walks on finite graphs}
\author{Vladislav Kargin}
\thanks{Department of Mathematics, Stanford University, CA 94305;
kargin@stanford.edu}
\date{June 2010}
\maketitle

\begin{center}
\textbf{Abstract}
\end{center}

\begin{quotation}
Several inequalities are proved for the mixing time of discrete-time quantum
walks on finite graphs. The mixing time is defined differently than in
Aharonov, Ambainis, Kempe and Vazirani (2001) and it is found that for
particular examples of walks on a cycle, a hypercube and a complete graph,
quantum walks provide no speed-up in mixing over the classical counterparts.
In addition, non-unitary quantum walks (i.e., walks with decoherence) are
considered and a criterion for their convergence to the unique stationary
distribution is derived.
\end{quotation}

\section{Introduction}

The origin of the concept of quantum walk lies in quantum computation
theory, where a quantum version of the classical random walk was invented in
an attempt to improve over classical computational algorithms. The early
papers that formulated the main ideas of quantum walk are \cite%
{aharonov_davidovich_zagury93} and \cite{meyer96}. Among numerous later
papers, we would like to point out \cite{farhi_gutmann98} where the
continuous-time quantum walk was defined and \cite%
{aharonov_ambainis_kempe_vazirani01} which defined and studied the
discrete-time quantum walk on finite graphs. An introductory review of
quantum walks written from the prospective of quantum computation can be
found in \cite{kempe03}. For recent developments the reader can also consult %
\cite{konno08}.

From the beginning, it became clear that quantum walks on both finite and
infinite graphs have many differences from the classical walk. For example,
the probability to find a particle at a particular vertex of a finite graph
does not converge to a limit but in general oscillates forever. However, the
average of this probability over time does converge to a limit, which can be
interpreted as follows. We start quantum walk in a certain state and measure
the particle at a random time $t,$ which is distributed uniformly over
interval $\left[ 0,T\right] .$ This measurement finds the particle at a
particular vertex $v$ with a probability $p\left( v,T\right) ,$ which
converges to a limit as $T\rightarrow \infty $. How large should $T$ be if
we want to make sure that $p\left( v,T\right) $ is close to its limit?

Let us introduce some definitions to make this question more precise.

The quantum walk on a finite graph is a 4-tuple $\left( G,S,\psi ,U\right) ,$
where $G=\left( V,E\right) $ is a finite graph, $S$ is a finite set, $\psi $
is a function in $L_{\mathbb{C}}^{2}\left( V\times S\right) $, and $U$ is a
unitary operator on $L^{2}\left( V\times S\right) .$ It is assumed that $%
\left\| \psi \right\| =1$. Elements of $S$ are called \emph{chiralities }and
the function $\psi _{t}=U^{t}\psi $ is the \emph{wave function }at time $%
t\in \mathbb{Z}$. If a measurement is performed over the system at time $t$,
then the walking particle is found at vertex $v$ in state $s$ with
probability $\left| \psi _{t}\left( v,s\right) \right| ^{2}.$

We assume that the quantum walk is \emph{local}$.$ That is, let $x$ and $%
x^{\prime }$ denote pairs $\left( v,s\right) $ and $\left( v^{\prime
},s^{\prime }\right) ,$ respectively. A quantum walk is local if $%
U_{x^{\prime }x}\equiv \left\langle \delta _{x^{\prime }},U\delta
_{x}\right\rangle \neq 0$ implies that $v\sim v^{\prime },$ that is,
vertices $v$ and $v^{\prime }$ are connected to each other. A local quantum
walk is called the \emph{general quantum walk} in \cite%
{aharonov_ambainis_kempe_vazirani01}.

A special case of the general quantum walk is the \emph{coined quantum walk }%
(\cite{aharonov_ambainis_kempe_vazirani01}). Here is how it is defined. Let $%
G$ be a $d$-regular graph and let $S=\left\{ 1,\ldots ,d\right\} .$ Assume
that the neighbors of each vertex $v$ are labelled as $v_{i}$ where $%
i=1,\ldots ,d.$ In addition, assume that if $v\neq w$ and $v_{i}=w_{j}$ then 
$i\neq j.$ (Such a labelling $L$ always exists on Cayley graphs of
finitely-generated groups, where we can identify elements of $S$ with
generators and inverses of generators of the group and write $v_{g}=vg$ and $%
v_{g^{-1}}=vg^{-1}$. In this case the choice of labelling is equivalent to
the choice of ordering of generators and their inverses.)

Define $U$ as follows. Let $x$ and $x^{\prime }$ denote pairs $\left(
v,s\right) $ and $\left( v^{\prime },s^{\prime }\right) ,$ respectively. If $%
v^{\prime }\nsim v,$ then $U_{x^{\prime }x}=0.$ Otherwise, $v^{\prime
}=v_{i} $ and $U_{x^{\prime }x}=\delta _{is}C_{s^{\prime }s},$ where $C$ is
a unitary matrix acting on $L^{2}\left( S\right) ,$ which is called the 
\emph{coin} of the quantum walk. It is easy to check that matrix $U$ is
unitary. Intuitively, let the particle be at vertex $v$ in state $\left|
s\right\rangle .$ Then, at the next moment the particle will be at vertex $%
v_{s}$ in the superposition state $C\left| s\right\rangle .$ This is the
coined quantum walk on $G$ corresponding to labelling $L$ and coin $C.$

A typical example of the coined quantum walk is the Hadamard quantum walk on
the cycle $\mathbb{Z}_{n}.$ In this case, the coin is the Hadamard
transformation: 
\begin{equation*}
C=\frac{1}{\sqrt{2}}\left( 
\begin{array}{cc}
1 & 1 \\ 
-1 & 1%
\end{array}%
\right) .
\end{equation*}%
Another popular choice of the coin is Grover's transformation: 
\begin{equation*}
C_{s^{\prime }s}=\left( 1-\frac{2}{d}\right) \delta _{s^{\prime }s}+\left( -%
\frac{2}{d}\right) \left( 1-\delta _{s^{\prime }s}\right) .
\end{equation*}%
That is, the state $s$ remains unchanged with amplitude $\left( 1-\frac{2}{d}%
\right) $ and moves to $s^{\prime }\neq s$ with amplitude $-2/d.$ We will
call walks with this coin the \emph{Grover quantum walks}.

A generalization of this concept is the \emph{non-unitary quantum walk }\cite%
{aharonov_ambainis_kempe_vazirani01}. A non-unitary quantum walk is
specified by $4$-tuple $\left( G,S,\rho ,\mathcal{T}\right) ,$ where $G$ and 
$S$ are as before a finite graph and a finite set, $\rho $ is a \emph{%
density matrix} (i.e., a positive unit-trace operator on $L_{\mathbb{C}%
}^{2}\left( V\times S\right) $), and $\mathcal{T}$ is a completely-positive
trace preserving operator acting on density matrices. In the literature, $%
\mathcal{T}$ is called a \emph{superoperator }\cite{preskill99}, or a \emph{%
quantum channel }\cite{preskill99}, or a trace-preserving \emph{quantum
operation }\cite{nielsen_chuang00}. We will use these terms as synonyms. Let 
$x$ denote a pair $\left( v,s\right) $. The probability to find a particle
at vertex $v$ in state $s$ at time $t$ is given by $\left\langle x|\mathcal{T%
}\rho |x\right\rangle .$ A non-unitary quantum walk is \emph{local} if $%
\left\langle x^{\prime }|\mathcal{T}\left( \rho \right) |x^{\prime
}\right\rangle >0$ for $x^{\prime }=\left( v^{\prime },s^{\prime }\right) $
implies that there is $x=\left( v,s\right) $ with $v\sim v^{\prime }$ such
that $\left\langle x|\rho |x\right\rangle >0.$ (The concept of locality is
more complicated in the non-unitary case and this definition is different
from the definition in \cite{aharonov_ambainis_kempe_vazirani01}.)

An example of a non-unitary quantum walk is given by a weighted sum of
unitary quantum walks. In this example, $\mathcal{T}\left( \rho \right)
=\sum_{i=1}^{k}p_{i}U_{i}\rho U_{i}^{\ast },$ where $U_{i}$ are unitary
operators, $p_{i}>0,$ and $\sum_{i=1}^{k}p_{i}=1.$ Intuitively, an operator $%
U_{i}$ is used at each step of the walk with probability $p_{i}.$ If all $%
U_{i}$ are local, then $\mathcal{T}$ is also local. Another example is $%
\mathcal{T}\left( \rho \right) =p\sum_{i=1}^{k}P_{i}\rho P_{i}+\left(
1-p\right) U\rho U^{\ast },$ where $P_{i}$ are projections and $%
\sum_{i=1}^{k}P_{i}=I.$ This is a walk in which with probability $p$ the
particle is measured and with probability $1-p$ it is evolved according to
the unitary operator $U$.

First, let us consider the case of unitary quantum walks. The probability
distribution $\left| \psi _{t}\left( x\right) \right| ^{2}$ in general does
not converge to any particular limit. Indeed, all eigenvalues of the matrix $%
U$ have unit absolute value. As a consequence, every eigenvector of $U$
corresponds to a stationary probability distribution. If the initial wave
function $\psi $ is a non-trivial superposition of the eigenvectors with
different eigenvalues, then $\psi _{t}$ continues to oscillate indefinitely.
In the classical case this phenomenon occurs only when the random walk
corresponds to a periodic Markov chain, and this case is not typical.

The time averages of the probabilities $\left| \psi _{t}\left( x\right)
\right| ^{2}$ do converge, and the limit 
\begin{equation*}
p\left( x\right) =\lim_{T\rightarrow \infty }\frac{1}{T}\sum_{t=0}^{T-1}%
\left| \psi _{t}\left( x\right) \right| ^{2}
\end{equation*}%
exists although may depend on the initial function $\psi $. We will call
this limit the \emph{time-averaged probability distribution} of the
particle. In order to quantify the convergence of the initial distribution
to this limit, let us define the distance of the initial distribution from
its time-average by the formula 
\begin{equation*}
d\left( T,\psi \right) =\sum_{x\in V\times S}\left| \frac{1}{T}%
\sum_{t=0}^{T-1}\left| \psi _{t}\left( x\right) \right| ^{2}-p\left(
x\right) \right| .
\end{equation*}%
This is the total variation distance between the averaged probability
distribution at time $T$ \ and its limit. By analogy with the classical
case, the \emph{mixing time} of a general quantum walk is defined as
follows: 
\begin{equation*}
t_{mix}\left( \varepsilon \right) =\sup_{\psi }\inf \left\{ T:d\left(
T^{\prime },\psi \right) \leq \varepsilon \text{ for all }T^{\prime }\geq
T\right\} .
\end{equation*}%
That is, this is the minimal time which is needed to reduce the distance
between the worst initial distribution and its time-averaged limit to a
quantity less than $\varepsilon .$

Another definition of the mixing time restricts the choice of initial wave
functions. Namely, 
\begin{equation*}
\widetilde{t_{mix}}\left( \varepsilon \right) =\sup_{\psi \in \mathcal{B}%
}\inf \left\{ T:d\left( T^{\prime }\right) \leq \varepsilon \text{ for all }%
T^{\prime }\geq T\right\} ,
\end{equation*}%
where $\mathcal{B}$ is the set of basis states, that is, $\psi \in \mathcal{B%
}$ if $\left| \psi \right| ^{2}$ is a delta-function concentrated at $\left(
v,s\right) $ \ This is the definition used in \cite%
{aharonov_ambainis_kempe_vazirani01}. Clearly, $\widetilde{t_{mix}}\left(
\varepsilon \right) \leq t_{mix}\left( \varepsilon \right) $. In the case of
classical random walks, these two mixing times are always equal to each
other. In the the case of quantum walks, they can be different.

We intend to estimate the mixing time in terms of the distance between
eigenvalues of $U.$ Let $\lambda _{k}=e^{i\beta _{k}},$ $k=1,\ldots ,m,$ be
the distinct eigenvalues of $U$. We define the distance between $\lambda _{k}
$ and $\lambda _{l}$ as the smallest distance along the unit circle: 
\begin{equation*}
d\left( \lambda _{k},\lambda _{l}\right) =\min \left\{ \left| \beta
_{k}-\beta _{l}+2\pi n\right| ,\text{ }n\in \mathbb{Z}\right\} .
\end{equation*}%
The \emph{relaxation time} of operator $U$ is defined as 
\begin{equation*}
t_{rel}=\max_{k\neq l}d\left( \lambda _{k},\lambda _{l}\right) ^{-1}.
\end{equation*}%
Finally, let us define the \emph{overlap} of two functions $\varphi $ and $%
\psi $ by the formula%
\begin{equation*}
Q\left( \varphi ,\psi \right) =\sum_{x\in V\times S}\left| \varphi \left(
x\right) \psi \left( x\right) \right| .
\end{equation*}

Note that if $\varphi $ and $\psi $ are two wave functions, then $0\leq
Q\leq 1$ by the Cauchy-Schwartz inequality.

\begin{theorem}
\label{theorem_disc_time_upper_bound} Let $U$ be the unitary transformation
on $L^{2}\left( V\times S\right) $ associated with a discrete-time quantum
walk (not necessarily coined). Let $U$ have $m$ distinct eigenvalues and the
relaxation time $t_{rel}$. Then,%
\begin{equation*}
t_{mix}\left( \varepsilon \right) \leq 2\pi \log \left( 2m\right) \frac{%
t_{rel}}{\varepsilon },
\end{equation*}
\end{theorem}

(The proofs of all theorems are in Appendix.)

It is interesting to compare this bound with the corresponding result for
the classical random walk, where $t_{mix}\left( \varepsilon \right) \leq
c\log \left( \frac{1}{\varepsilon \pi _{\min }}\right) t_{rel},$ where $\pi
_{\min }$ is the smallest probability in the limit distribution (see for
example Theorems 12.3 and 12.4 on p. 155 in \cite{levin_peres_wilmer09}). In
many cases the limit distribution is uniform and this bound can be written
as $t_{mix}\left( \varepsilon \right) \leq c\log \left( n/\varepsilon
\right) t_{rel},$ where $n$ is the number of vertices in the graph. Note,
however, that $t_{rel}$ have a different meaning in the classical case where
it denotes the inverse of the difference between $1$ (the largest
eigenvalue) and the second largest eigenvalue (i.e., the inverse of the
``spectral gap'').

Another significant difference in the formulas for the mixing time is that $%
\varepsilon $ enters as $\log \varepsilon $ and $\varepsilon ^{-1}$ in the
classical and quantum cases, respectively. This is due to the fact that the
convergence is exponentially fast in the classical case and polynomial (even
linear) in the quantum case.

Finally, it is worthwhile to note that in many cases the classical bound $%
t_{mix}\leq c\log \left( n\right) t_{rel}$ is not optimal, and a large
literature is devoted to improvement of this result to $t_{mix}\leq ct_{rel}$
with a sharp constant $c.$

For the lower bound we prove the following result.

\begin{theorem}
\label{theorem_disc_time_lower_bound} Let $U$ be the unitary transformation
on $L^{2}\left( V\times S\right) $ associated with a discrete-time quantum
walk (not necessarily coined). Suppose that $U$ has only real eigenvectors.
Let $\lambda $ and $\lambda ^{\prime }$ be two distinct eigenvalues with the
corresponding eigenvectors $\psi $ and $\psi ^{\prime }.$ Assume that $%
d\left( \lambda ,\lambda ^{\prime }\right) \leq 2$ and $\varepsilon \leq
Q\left( \psi ,\psi ^{\prime }\right) /80.$ Then,%
\begin{equation*}
t_{mix}\left( \varepsilon \right) \geq \frac{Q\left( \psi ,\psi ^{\prime
}\right) }{8\varepsilon d\left( \lambda ,\lambda ^{\prime }\right) }.
\end{equation*}
\end{theorem}

In particular, if $\lambda \,$\ and $\lambda ^{\prime }$ are two eigenvalues
with the smallest distance between them along the circle, then $d\left(
\lambda ,\lambda ^{\prime }\right) ^{-1}=t_{rel}$ and we obtain the estimate 
\begin{equation*}
t_{mix}\left( \varepsilon \right) \geq \frac{1}{8}Q\left( \psi ,\psi
^{\prime }\right) \frac{t_{rel}}{\varepsilon }
\end{equation*}%
valid for all sufficiently small $\varepsilon .$

The main message of Theorems \ref{theorem_disc_time_upper_bound} and \ref%
{theorem_disc_time_lower_bound} is that the relation of the mixing and
relaxation times in the quantum case is similar to the analogous relation in
the classical case. However, the relaxation time is defined differently in
the quantum case. It is not the inverse of the difference between the
largest and the second largest eigenvalue, but the inverse of the minimal
distance between all distinct eigenvalues.

Previously, the speed of convergence of (unitary) discrete-time quantum
walks was investigated in \cite{aharonov_ambainis_kempe_vazirani01}. The
upper bound for the quantum walks that we obtain in Theorem \ref%
{theorem_disc_time_upper_bound} is similar to the bound in Theorem 6.1 of %
\cite{aharonov_ambainis_kempe_vazirani01}. The mixing time is $O\left(
t_{rel}\log \left( m\right) \right) $ where $t_{rel}$ is the inverse of the
minimal distance between the distinct eigenvalues of the matrix $U$. The
main difference of our result from the result in \cite%
{aharonov_ambainis_kempe_vazirani01} is that we have $\log \left( m\right) $
instead of $\log (n),$ where $m$ and $n$ are the numbers of distinct and all
eigenvalues, respectively.

This difference is significant for the case of the discrete walk on the
hypercube, where the number of eigenvalues is $2^{d}$ and the number of
distinct eigenvalues is $d+1.$ In particular, we show that the mixing time
on the hypercube is $O(n\log n/\varepsilon )$ and not exponential as was
suggested in \cite{moore_russell02} based on previous estimates in \cite%
{aharonov_ambainis_kempe_vazirani01}.

The lower bound that we obtain is in terms of the relaxation time $t_{rel}$.
It essentially says that the mixing time is $\Omega \left( t_{rel}\right) .$
This bound is different from the bound obtained in \cite%
{aharonov_ambainis_kempe_vazirani01}, which is formulated in terms of a
geometrical property of the underlying graph. In addition, the mixing time
is defined differently in \cite{aharonov_ambainis_kempe_vazirani01}. As a
result, the mixing time of the Hadamard walk on the cycle is of order $%
O\left( n\log n/\varepsilon ^{3}\right) $ in Theorem 4.2 of \cite%
{aharonov_ambainis_kempe_vazirani01}, and of order $O\left( n^{2}\log
n/\varepsilon \right) $ in our Example 1. In the classical case, the mixing
time is of the order $O\left( n^{2}\log \left( \varepsilon ^{-1}\right)
\right) .$

Now let us consider non-unitary quantum walks. The study of these walks
helps us to understand how the decoherence affects performance of quantum
algorithms. It was noted (see \cite{kendon_tregenna03}) that decoherence in
quantum walks can be useful for quantum algorithms. In particular, it
appears that a small amount of decoherence can speed up the mixing of the
walk. Numeric evidence in \cite{kendon_tregenna03} was later corroborated by
analytical estimates in \cite{richter07}. More information about decoherence
in quantum walks and additional references can be found in the review
article \cite{kendon07}.

Let $\mathcal{M}$ denote the linear space of Hermitian linear operators
acting on $L_{\mathbb{C}}^{2}\left( V\times S\right) .$ The space $\mathcal{M%
}$ is a Hilbert space with respect to the norm $\left\| \rho \right\| _{2}=%
\left[ \mathrm{Tr}\left( \rho ^{2}\right) \right] ^{1/2}$ (which we call $%
L^{2}$\emph{-norm}). Other useful norms on $\mathcal{M}$ are $\left\| \rho
\right\| _{1}=\mathrm{Tr}\left( \left| \rho \right| \right) $ where $\left|
\rho \right| =\sqrt{\rho ^{2}}$ and $\left\| \rho \right\| _{L^{2}\left( \mu
\right) }=\left[ \mathrm{Tr}\left( \mu \rho ^{2}\right) \right] ^{1/2},$
where $\mu $ is a density matrix. We call these norms the \emph{trace} and $%
L^{2}\left( \mu \right) $ norms, respectively. Superoperators are operators
on $\mathcal{M}$ which possess some additional properties. Some well-known
properties of superoperators are summarized in the proposition below.

\begin{proposition}
\label{proposition_properties_channels}Superoperator $\mathcal{T}$ is a
contraction in the trace norm (i.e., $\left\| \mathcal{T}\left( \rho \right)
\right\| _{1}\leq \left\| \rho \right\| _{1}$). There exists a density
matrix $\rho _{st}$ such that $\mathcal{T}\rho _{st}=\rho _{st}.$
\end{proposition}

This proposition is an immediate consequence of Theorem 9.2 and Exercise 9.9
in \cite{nielsen_chuang00}.

Note that in many cases $\mathcal{T}$ is not self-adjoint in $L^{2}$ norm.
Moreover, recall that in the classical case the stochastic matrix of a
random walk is always self-adjoint with respect to the norm $L^{2}\left( \mu
\right) ,$ where $\mu $ is the stationary probability distribution. (This
result can be traced to the fact that every random walk is a \emph{reversible%
} Markov chain.) In contrast, the superoperator of a non-unitary quantum
walk $\mathcal{T}$ is not necessarily self-adjoint with respect to the norm $%
L^{2}\left( \rho _{st}\right) $. In fact, it appears that $\mathcal{T}$ is
not even a normal operator (i.e., $\mathcal{T}^{\ast }\mathcal{T\neq TT}%
^{\ast })$ in many situations of interest.

Proposition \ref{proposition_properties_channels} establishes the existence
of the stationary density matrix. However, it does not say anything about
the uniqueness or convergence properties, and we cannot expect that these
properties hold in general. For example, a unitary quantum walk typically
has many stationary density matrices and the convergence fails unless we
average density matrices over time. The following theorem establishes the
uniqueness and convergence properties provided that the quantum walk
satisfies a certain condition. Let us call a density matrix $\rho $ \emph{%
strictly positive} and write $\rho >0,$ if $\left\langle x,\rho
x\right\rangle =0$ implies that $x=0.$ Next, let $\mathcal{T}$ be a linear
operator acting on $\mathcal{M}$. We will call $\mathcal{T}$ \emph{strongly
positive} if for every density matrix $\rho $ there exists an integer $n>0$
such that $\mathcal{T}^{n}\rho >0$.

(This definition is similar to a corresponding definition in the theory of
Markov chains, in which it is shown that a stochastic matrix of a Markov
chain is strongly positive if and only if the Markov chain is ergodic, that
is, aperiodic and irreducible.)

The \emph{multiplicity} of an eigenvalue $\lambda $ is defined as $\dim \ker
\left( \lambda I-T\right) .$ The \emph{rank} of $\lambda $ is $%
\sup_{p>0}\dim \ker \left( \lambda I-T\right) ^{p}.$ The eigenvalue is
called simple if its rank equals $1.$

\begin{theorem}
\label{theorem_nonunitary_convergence}Let $\mathcal{T}$ be a strongly
positive superoperator. Then, (i) $\mathcal{T}$ has a simple eigenvalue $1$.
(ii) The corresponding eigenvector $\rho _{st}$ is a strictly positive
density matrix. (iii) For every initial density matrix $\rho _{0},$ $%
\mathcal{T}^{n}\rho _{0}\rightarrow \rho _{st},$ as $n\rightarrow \infty .$
\end{theorem}

Proof is in Appendix.

After the convergence to the stationary distribution is established, it is
natural to ask for an estimate on the mixing time. First, let us define the
mixing time for a non-unitary quantum walks. The definition is different
from the definition for the unitary walks since no time-averaging is
necessary. The measurement at time $t$ finds the walking particle at the
vertex $v$ in state $s$ with probability $p_{t}\left( x,\rho \right)
=\left\langle x|\mathcal{T}^{t}\left( \rho \right) |x\right\rangle $, where $%
x$ denote the pair $\left( v,s\right) $ and $\rho $ is the initial density
matrix$.$ If $\mathcal{T}$ is strongly positive, then these probabilities
converge to a limit $p\left( x\right) =\left\langle x|\rho
_{st}|x\right\rangle ,$ which does not depend on the initial density matrix.
Hence, we can define the total variation distance as $d\left( t,\rho \right)
=\sum_{x\in V\times S}\left| \left\langle x|\mathcal{T}^{t}\left( \rho
\right) |x\right\rangle -p\left( x\right) \right| $. The corresponding
mixing time can be defined as 
\begin{equation*}
t_{mix}\left( \varepsilon \right) =\sup_{\rho }\inf \left\{ t:d\left(
t^{\prime },\rho \right) \leq \varepsilon \text{ for all }t^{\prime }\geq
t\right\} .
\end{equation*}

Unfortunately, while it is easy to see that the asymptotic behavior of $%
\mathcal{T}^{t}$ is governed by the spectral radius of $\mathcal{T},$ it is
difficult to estimate the mixing time because of the non-normality of
operator $\mathcal{T}$. The essential difficulty is that for such operators
it is hard to estimate the duration of the transient behavior. It is the
same problem that makes it difficult to estimate the mixing time for
non-reversible Markov chains.

(In one particular example of a non-unitary continuous-time walk on cycle
this difficulty has been overcome and an estimate on the mixing time has
been derived in \cite{richter07}.)

We consider several examples of unitary walks in this paper. The table
summarizes results for unitary quantum walks on a complete graph, a cycle,
and a hypercube. 
\begin{equation*}
\begin{tabular}{ll}
& Mixing time \\ 
Complete graph with $n$ vertices & $\frac{c_{1}}{\varepsilon }\leq
t_{mix}\left( \varepsilon \right) \leq \frac{c_{2}}{\varepsilon }$ \\ 
Cycle with $n$ vertices & $\frac{c_{1}}{\varepsilon .}n^{2}\leq
t_{mix}\left( \varepsilon \right) \leq \frac{c_{2}}{\varepsilon }n^{2}\log n$
\\ 
Hypercube with $2^{n}$ vertices & $\frac{1}{2\varepsilon }n\leq
t_{mix}\left( \varepsilon \right) \leq \frac{2\pi }{\varepsilon }n\log n$%
\end{tabular}%
\end{equation*}

It appears from this table that the mixing time for quantum walks is of
similar order as that for the corresponding classical random walks. In
particular, the unitary quantum walks do not allow a quadratic speedup over
classical walks, in contrast to the results for the mixing time in \cite%
{aharonov_ambainis_kempe_vazirani01}. The reason for this difference is that
the mixing time defined in \cite{aharonov_ambainis_kempe_vazirani01}
restricts the initial distributions of the particle to the class of
distributions concentrated on a particular vertex of a graph, while we allow
for arbitrary initial distributions. Note that this result does not rule out
that the quadratic speedup can be achieved by non-unitary quantum walks.
Some evidence in favour of this conjecture can be found in \cite%
{kendon_tregenna03} and \cite{richter07}.

The rest of the paper is organized as follows. In the next section, we apply
bounds on mixing times to particular examples of quantum walks on the cycle,
hypercube, and complete graph. The proofs of the theorems are relegated to
Appendix.

\section{Examples}

\textbf{Example 1}. (\textbf{Cycle}) .

\begin{proposition}
The mixing time for the Hadamard quantum walk on the $n$-cycle satisfies the
following inequalities: 
\begin{equation*}
\frac{c_{1}}{\varepsilon }n^{2}\leq t_{mix}\left( \varepsilon \right) \leq 
\frac{c_{2}}{\varepsilon }n^{2}\log \left( n\right) ,
\end{equation*}%
where $c_{1}$ and $c_{2}$ are positive constants.
\end{proposition}

\textbf{Proof: }The eigenvalues of the Hadamard walk on the cycle with $n$
vertices were found in \cite{aharonov_ambainis_kempe_vazirani01}. They are 
\begin{equation*}
t_{k}^{\left( 1,2\right) }=\frac{1}{\sqrt{2}}\left( \cos \left( \frac{2\pi }{%
n}k\right) \pm i\sqrt{1+\sin ^{2}\left( \frac{2\pi }{n}k\right) }\right) ,
\end{equation*}%
where $k=0,\ldots ,n-1.$ In order to describe the eigenvectors, let $\chi
_{k}$, $0\leq k\leq n-1,$ be functions in $L^{2}\left( \mathbb{Z}_{n}\right) 
$ defined by the formula $\chi _{k}=\sum_{r=0}^{n-1}\exp \left( 2\pi i\frac{%
kr}{n}\right) \delta _{r}.$ Then all eigenvectors have the form $v\otimes
\chi _{k}$ where $v$ is a $2$-vector that depends on $k.$

Indeed, if $S$ and $S^{\ast }$ are the left and right shift operator on $%
L^{2}\left( \mathbb{Z}_{n}\right) $, respectively, then we can write $U$ as
a $2$-by-$2$ block matrix, with blocks $U_{11}$ and $U_{12}$ equal $S/\sqrt{2%
},$ and blocks $U_{21}$ and $U_{22}$ equal $-S^{\ast }/\sqrt{2}$ and $%
S^{\ast }/\sqrt{2},$ respectively. It follows that $U\left( v\otimes \chi
_{k}\right) =A_{k}v\otimes \chi _{k},$ where 
\begin{equation*}
A_{k}=\frac{1}{\sqrt{2}}\left( 
\begin{array}{cc}
\exp \left( -2\pi i\frac{k}{n}\right)  & \exp \left( -2\pi i\frac{k}{n}%
\right)  \\ 
-\exp \left( 2\pi i\frac{k}{n}\right)  & \exp \left( 2\pi i\frac{k}{n}%
\right) 
\end{array}%
\right) ,
\end{equation*}%
Let 
\begin{equation*}
t_{k}^{\left( 1,2\right) }=\frac{1}{\sqrt{2}}\left( \cos \left( \frac{2\pi }{%
n}k\right) \pm i\sqrt{1+\sin ^{2}\left( \frac{2\pi }{n}k\right) }\right) .
\end{equation*}%
Then, eigenvectors of $A_{k}$ can be written as $v_{k}^{\left( \alpha
\right) }=\left( c,1\right) $ with $c=1-\sqrt{2}t_{k}^{\left( \alpha \right)
}\exp \left( -2\pi ik/n\right) .$The corresponding eigenvectors of $U$ are $%
v_{k}^{\left( \alpha \right) }\otimes \chi _{k}$ with eigenvalues $\lambda
_{k}^{\left( \alpha \right) }=t_{k}^{\left( \alpha \right) }$ for $\alpha
=1,2.$

Note that%
\begin{equation*}
\beta _{k}^{\left( 1,2\right) }:=\arg \left( t_{k}^{\left( 1,2\right)
}\right) =\pm \arccos \left[ \frac{1}{\sqrt{2}}\cos \left( \frac{2\pi }{n}%
k\right) \right] 
\end{equation*}%
The smallest difference between $\beta _{k}$ occurs when $k=0$ and $1$ and
it can be estimated by $c/n^{2}$ for a suitable constant $c.$ It follows
that the relaxation time is $t_{rel}\sim cn^{2},$ and by Theorem \ref%
{theorem_disc_time_upper_bound} the mixing time is 
\begin{equation*}
t_{mix}\left( \varepsilon \right) \leq \frac{c_{2}}{\varepsilon }n^{2}\log
\left( n\right) ,
\end{equation*}%
with a certain constant $c_{2}>0.$

It is easy to estimate the overlap of eigenvectors that correspond to
eigenvalues $\beta _{0}$ and $\beta _{1}.$ It is greater than $0.97$ for all 
$n.$ By Theorem \ref{theorem_disc_time_lower_bound}, we have 
\begin{equation*}
\frac{c_{1}}{\varepsilon }n^{2}\leq t_{mix}\left( \varepsilon \right) .
\end{equation*}%
QED.

\textbf{Example 2. (Hypercube)}

The mixing time of the quantum walk on a hypercube was previously studied in %
\cite{moore_russell02}, and we use their setup in the definition of quantum
walk. The quantum walk on the hypercube is also analyzed in \cite{kempe05}
with emphasis on the hitting time of the walk.

Consider a hypercube graph $\left( \mathbb{Z}_{2}\right) ^{n}$ with $2^{n}$
vertices. We think about vertices as indexed by numbers from $0$ to $2^{n}-1$
in the binary representation with $n$ digits. The edges of the graph are put
between numbers that are different in one bit only. The set of states $S$
consists of $n$ elements$.$ We consider the Grover quantum walk. That is, a
particle at vertex $v$ in state $s$ goes to the vertex $w$ which is
different from vertex $v$ only in the bit $s.$ It remains in state $s$ with
amplitude $2/n-1$ and goes to state $s^{\prime }$ with amplitude $2/n.$

\begin{proposition}
The mixing time for the Grover quantum walk on the $n$-dimensional hypercube
satisfies the following inequalities: 
\begin{equation*}
\frac{1}{2\varepsilon }n\leq t_{mix}\left( \varepsilon \right) \leq \frac{%
2\pi }{\varepsilon }n\log \left( n\right) .
\end{equation*}
\end{proposition}

\textbf{Proof: }The eigenvalues and eigenvectors of the Grover quantum walk
on the $n$-dimensional hypercube were found by Moore and Russell in \cite%
{moore_russell02}. The eigenvalues are 
\begin{equation*}
\lambda _{k}^{\pm }=1-\frac{2k}{n}\pm 2i\frac{\sqrt{k(n-k)}}{n},
\end{equation*}%
where $k=0,\ldots ,n.$ We describe eigenvectors below. For the convenience
of the reader, we also give a short verification of the result .

For each sequence $t=\left( t_{1},t_{2},\ldots ,t_{n}\right) $ of $0$ and $%
1, $ define $\chi _{t}\in L^{2}\left( \mathbb{Z}_{2}^{n}\right) $ by the
formula $\chi _{t}\left( x_{1},\ldots ,x_{n}\right) =2^{-n/2}\left(
-1\right) ^{\sum t_{i}x_{i}}.$ All eigenvectors of the matrix $U$ have the
form $v_{t}^{\left( j\right) }\otimes \chi _{t},$ where $v_{t}^{\left(
j\right) }$ is an $n$-vector that depends on $t,$ and $j=1,...,n.$

Indeed, the unitary matrix $U$ can be written as a $n$-by-$n$ block matrix,
in which the $ij$-th block is $bS_{j}$ if $i\neq j$ and $aS_{j}$ if $i=j.$
Here $a=2/n-1,$ $b=2/n$ and $S_{k}:L^{2}\left( \mathbb{Z}_{2}^{n}\right)
\rightarrow L^{2}\left( \mathbb{Z}_{2}^{n}\right) $ is the shift operator
which acts as follows: $\left( S_{k}f\right) \left( x_{1},\ldots
,x_{n}\right) =f\left( x_{1},\ldots ,x_{k}+1,\ldots ,x_{n}\right) ,$ where
addition is modulo $2.$ Note that $S_{k}\chi _{t}=\left( -1\right)
^{t_{k}}\chi _{t}.$

A computations shows that $U\left( v\otimes \chi _{t}\right) =A\left(
v\right) \otimes \chi _{t},$ where $A$ is an $n$-by-$n$ matrix (which
depends on $t$) with entries $A_{ij}=\left[ \delta _{ij}a+\left( 1-\delta
_{ij}\right) b\right] \left( -1\right) ^{t_{j}}.$ In other form, $A=D+bP,$
where $D=-\mathrm{diag}\left( \left( -1\right) ^{t_{1}},\left( -1\right)
^{t_{2}},\ldots ,\left( -1\right) ^{t_{n}}\right) ,$ and $P=\left|
1,1,\ldots ,1\right\rangle \left\langle \left( -1\right) ^{t_{1}},\left(
-1\right) ^{t_{2}},\ldots ,\left( -1\right) ^{t_{n}}\right| .$

It is easy to verify that the following vectors are eigenvectors of $A.$ Let 
$k$ be the number of non-zero entries in vector $t$. First, assume that $%
n>k\geq 1$ and define $x=\pm i\sqrt{\frac{k}{n-k}}$ and $v_{r}=x^{1-t_{r}}.$
Then $v=\left( v_{1},\ldots ,\nu _{n}\right) $ is an eigenvector of $A$ with
eigenvalue 
\begin{equation*}
\lambda _{k}^{\pm }=1-\frac{2k}{n}\pm 2i\frac{\sqrt{k(n-k)}}{n}.
\end{equation*}

In addition, note that every non-zero vector $v$ such that $v_{r}=0$ if $%
t_{r}=0$ and $\sum_{r=1}^{n}v_{r}=0$, is an eigenvector of $A$ with
eigenvalue $1.$ The set of such vectors form an eigenspace of dimension $k-1$%
. Similarly, every non-zero $v$ such that $v_{r}=0$ if $t_{r}=1$ and $%
\sum_{r=1}^{n}v_{r}=0$ is an eigenvector of $A$ with eigenvalue $-1.$ The
set of such vectors forms an eigenspace of dimension $n-k-1.$

For the case when $t=0,$ the vector $\left( 1,\ldots ,1\right) $ is an
eigenvector with eigenvalue $1$ and its orthogonal complement is eigenspace
of $-1.$ For $t=\left( 1,1,...,1\right) ,$ the situation is reverse.

By counting dimensions of eigenspaces, it is clear that these are all
eigenvalues of matrix $A.$ Since there are $2^{n}$ different choices of
vector $t,$ we also found all eigenvalues of matrix $U.$ It follows that
these eigenvalues are $\pm 1,$ and $\lambda _{k}^{\pm }$ for $k=1,\ldots
,n-1.$

From the formula for eigenvalues, the distance between distinct eigenvalues
can be estimated from below as $\Delta >\frac{2}{n}.$ Hence, $t_{rel}<\frac{n%
}{2}.$

By applying Theorem \ref{theorem_disc_time_upper_bound}, we find 
\begin{equation*}
t_{mix}\left( \varepsilon \right) \leq \frac{2\pi }{\varepsilon }n\log
\left( n\right) .
\end{equation*}

For the lower bound, consider for simplicity the case of even $n=2m.$ (The
case of odd $n$ is similar.) Let $v_{x,k}$ denote the value of function $%
v\in L^{2}\left( \left( \mathbb{Z}_{2}\right) ^{n}\times \mathbb{Z}%
_{n}\right) $ on vertex $x$ and state $k,$ and consider the eigenvectors
that correspond to eigenvalues $\lambda _{m}^{+}$ and $\lambda _{m+1}^{+}$
respectively. Then, it is easy to compute the overlap of these eigenvectors
as $\sqrt{1-\left( m+1\right) /\left( 2m^{3}\right) }\sim 1$ for large $n.$
The distance between arguments of eigenvalues $\lambda _{m}^{+}$ and $%
\lambda _{m+1}^{-}$ is approximately $2/n.$ Hence, by Theorem \ref%
{theorem_disc_time_lower_bound} we have the inequality%
\begin{equation*}
t_{mix}\left( \varepsilon \right) \gtrsim \frac{n}{2\varepsilon }.
\end{equation*}%
QED.

\textbf{Example 3. (Complete graph)}

There are several ways to define a discrete-time walk on the complete graph
with $n$ vertices. We will consider the following variant. Let $\left|
S\right| =n.$ Define the entries of the unitary matrix as follows. 
\begin{equation*}
U\left( ws^{\prime },vs\right) =\delta _{ws}\left\{ \left( 1-\frac{2}{n}%
\right) \delta _{vs^{\prime }}+\left( -\frac{2}{n}\right) \left( 1-\delta
_{vs^{\prime }}\right) \right\} .
\end{equation*}%
In words, let the particle start at vertex $v$ in state $s.$ Then at the
next moment of time it will be at vertex $w=s.$ The particle moves to state $%
s^{\prime }$ with amplitude $\left( -2/n\right) $ for $s^{\prime }\neq v.$
If $s^{\prime }=v,$ then the amplitude of the transition is $\left(
1-2/n\right) .$

\begin{proposition}
The mixing time for the quantum walk on the complete graph satisfies
inequalities:%
\begin{equation*}
\frac{c_{1}}{\varepsilon }\leq t_{mix}\left( \varepsilon \right) \leq \frac{%
c_{2}}{\varepsilon },
\end{equation*}%
where $c$ is a positive constant.
\end{proposition}

\textbf{Proof:} We will show that the eigenvalues of $U$ are $1,$ $-1,$ $i$
and $-i$ with multiplicities $n(n-1)/2,$ $1+(n-1)(n-2)/2,$ $n-1,$ and $%
n-1,\, $respectively. Let $X_{vs}=\psi \left( v,s\right) .$ Then, the action
of $U$ can be written as follows: 
\begin{equation*}
U:X\rightarrow X^{T}\left( I-\frac{2}{n}1_{n}1_{n}^{T}\right) ,
\end{equation*}%
where $X^{T}$ is the transposed matrix $X,$ $1_{n}$ is the column $n$-by-$1$
vector that consists of all ones, and $1_{n}^{T}$ is the corresponding row
vector.

If $X=X^{T}$ and all columns of $X$ sum to $0,$ then $U(X)=X.$ This gives us
an eigenspace of operator $U$ with eigenvalue $1$ and dimension $n(n-1)/2.$

Similarly if $X^{T}=-X$ and all columns of $X$ sum to $0,$ then $U\left(
X\right) =-X.$ This gives us an eigenspace of $U$ with eigenvalue $-1.$ In
addition, $U\left( I\right) =-I.$ Hence, the dimension of the eigenspace
with eigenvalue $-1$ is $1+(n-1)(n-2)/2$. In order to find the eigenspaces
with eigenvalues $\pm i,$ consider $U^{2}.$ It acts as follows:%
\begin{equation*}
U^{2}\text{: }X\rightarrow \left( I-\frac{2}{n}1_{n}1_{n}^{T}\right) X\left(
I-\frac{2}{n}1_{n}1_{n}^{T}\right) .
\end{equation*}%
Let $c_{1},\ldots ,c_{n}$ and $r_{1},\ldots ,r_{n}$ be arbitrary numbers
satisfying the conditions $\sum c_{i}=\sum r_{j}=0,$ and define $%
X_{ij}=r_{i}+c_{j}.$ Then $U^{2}\left( X\right) =-X.$ Hence, these matrices
belong to the eigenspace of $U^{2}$ with eigenvalue $-1.$ The dimension of
this space is $2n-2.$ It follows that $U$ has two eigenspaces of dimension $%
n-1$ which correspond to eigenvalues $i$ and $-i,$ respectively. By counting
dimensions we confirm that we have found all eigenvalues and eigenspaces of
matrix $U.$

The relaxation time is $t_{rel}=1/\sqrt{2}.$ By applying Theorems \ref%
{theorem_disc_time_lower_bound} and \ref{theorem_disc_time_upper_bound}, the
mixing time satisfies inequalities

\begin{equation*}
\frac{c_{1}}{\varepsilon }\leq t_{mix}\left( \varepsilon \right) \leq \frac{%
c_{2}}{\varepsilon }
\end{equation*}%
for some positive $c_{1}$ and $c_{2}.$

\appendix

\section{Proofs of Theorems}

\textbf{Proof of Theorem \ref{theorem_disc_time_upper_bound}:} Let $V_{k}$
be the eigenspace corresponding to eigenvalue $\beta _{k}$ of operator $U$.
Then, we can write 
\begin{equation*}
\psi =\sum_{k=1}^{m}c_{k}\varphi _{k},
\end{equation*}%
where $\varphi _{k}\in V_{k}$, $\left\| \varphi _{k}\right\| =1,$ $%
\sum_{k}\left| c_{k}\right| ^{2}=1.$ Then at time $t,$ 
\begin{equation*}
\psi \left( x,t\right) =\sum_{k=1}^{m}c_{k}\varphi _{k}\left( x\right)
\left( \lambda _{k}\right) ^{t},
\end{equation*}%
and 
\begin{equation*}
\left| \psi \left( x,t\right) \right| ^{2}=\sum_{k=1}^{m}\left| c_{k}\right|
^{2}\left| \varphi _{k}\left( x\right) \right| ^{2}+\sum_{k\neq l}c_{k}%
\overline{c_{l}}\varphi _{k}\left( x\right) \overline{\varphi }_{l}\left(
x\right) \left( \lambda _{k}\overline{\lambda }_{l}\right) ^{t}.
\end{equation*}%
Hence, 
\begin{eqnarray*}
\frac{1}{T}\sum_{t=0}^{T-1}\left| \psi \left( x,t\right) \right| ^{2}dt
&=&\sum_{k=1}^{m}\left| c_{k}\right| ^{2}\left| \varphi _{k}\left( x\right)
\right| ^{2} \\
&&+\frac{1}{T}\sum_{k\neq l}c_{k}\overline{c_{l}}\varphi _{k}\left( x\right) 
\overline{\varphi }_{l}\left( x\right) \frac{\left( \lambda _{k}\overline{%
\lambda }_{l}\right) ^{T}-1}{\lambda _{k}\overline{\lambda }_{l}-1}.
\end{eqnarray*}

It follows that 
\begin{equation*}
p\left( x\right) :=\lim_{T\rightarrow \infty }\frac{1}{T}\sum_{t=1}^{T-1}%
\left| \psi \left( x,t\right) \right| ^{2}=\sum_{k=1}^{m}\left| c_{k}\right|
^{2}\left| \varphi _{k}\left( x\right) \right| ^{2},
\end{equation*}%
and 
\begin{equation*}
d\left( T\right) =\frac{1}{T}\sum_{x\in G}\left| \sum_{k\neq l}c_{k}%
\overline{c_{l}}\varphi _{k}\left( x\right) \overline{\varphi }_{l}\left(
x\right) \frac{\left( \lambda _{k}\overline{\lambda }_{l}\right) ^{T}-1}{%
\lambda _{k}\overline{\lambda }_{l}-1}\right| .
\end{equation*}

In order to bound this quantity, note that 
\begin{equation*}
\sum_{x\in G}\left| c_{k}\overline{c_{l}}\varphi _{k}\left( x\right) 
\overline{\varphi }_{l}\left( x\right) \right| \leq Q\left( \varphi
_{k},\varphi _{l}\right) \left| c_{k}\right| \left| c_{l}\right| \leq \left|
c_{k}\right| \left| c_{l}\right| ,
\end{equation*}%
and, therefore, 
\begin{equation*}
d\left( T\right) \leq \frac{2}{T}\sum_{k\neq l}\frac{\left| c_{k}\right|
\left| c_{l}\right| }{\left| e^{i\left( \beta _{k}-\beta _{l}\right)
}-1\right| }.
\end{equation*}%
Note that $\left| e^{i\left( \beta _{k}-\beta _{l}\right) }-1\right| \geq 
\frac{2}{\pi }d\left( \beta _{k},\beta _{l}\right) $ where $d\left( \beta
_{k},\beta _{l}\right) $ is the distance between $\beta _{k}$ and $\beta _{l}
$ modulo $2\pi ,$ that is, $d\left( \beta _{k},\beta _{l}\right)
:=\min_{s}\left\{ \left| \beta _{l}-\beta _{k}+2\pi s\right| \right\} .$ Let 
$\Delta $ denote $\min_{k,l}d\left( \beta _{k},\beta _{l}\right) $ and
assume that $0\leq \beta _{1}<\beta _{2}<\ldots <\beta _{m}<2\pi .$ Then we
can write $d\left( \beta _{k},\beta _{l}\right) \geq \Delta d\left(
k,l\right) ,$ where $d\left( k,l\right) =\min \left\{ \left| k-l\right|
,\left| k-l+m\right| ,\left| k-l-m\right| \right\} .$ This inequality holds
because the shortest arc of the circle between $\beta _{k}$ and $\beta _{l}$
contains $d\left( k,l\right) $ non-overlapping intervals whose endpoints are 
$\beta _{s}$ and the length of each of these intervals is at least $\Delta .$
It follows that $\left| e^{i\left( \beta _{k}-\beta _{l}\right) }-1\right|
\geq \frac{2\Delta }{\pi }d\left( k,l\right) .$

Next,  
\begin{eqnarray*}
d\left( T\right)  &\leq &\frac{\pi }{T\Delta }\sum_{k\neq l}\frac{\left|
c_{k}\right| \left| c_{l}\right| }{d\left( k,l\right) } \\
&\leq &\frac{\pi }{T\Delta }\sum_{k\neq l}\frac{\frac{1}{2}\left( \left|
c_{k}\right| ^{2}+\left| c_{l}\right| ^{2}\right) }{d\left( k,l\right) }=%
\frac{\pi }{T\Delta }\sum_{k\neq l}\frac{\left| c_{k}\right| ^{2}}{d\left(
k,l\right) }.
\end{eqnarray*}%
This sum can be estimated as follows: 
\begin{equation*}
\sum_{k\neq l}\frac{\left| c_{k}\right| ^{2}}{d\left( k,l\right) }\leq
\sum_{k=1}^{m}\left| c_{k}\right| ^{2}\left( \sum_{l\neq k}\frac{1}{d\left(
k,l\right) }\right) =\sum_{k=1}^{m}\left| c_{k}\right| ^{2}\left(
\sum_{l=2}^{m}\frac{1}{d\left( 1,l\right) }\right) .
\end{equation*}%
We estimate 
\begin{equation*}
\sum_{l=2}^{m}\frac{1}{d\left( 1,l\right) }\leq 2\sum_{k=1}^{\left[ m/2%
\right] }\frac{1}{k}\leq 2\left( 1+\log \left[ \frac{m}{2}\right] \right)
<2\log \left( 2m\right) 
\end{equation*}%
and since $\sum_{k=1}^{m}\left| c_{k}\right| ^{2}=1,$ we have 
\begin{equation*}
d\left( T\right) \leq \frac{2\pi }{T\Delta }\log \left( 2m\right) 
\end{equation*}

Hence, 
\begin{equation*}
t_{mix}\left( \varepsilon \right) \leq 2\pi \log \left( 2m\right) \frac{%
t_{rel}}{\varepsilon }.
\end{equation*}%
QED.

\textbf{Proof of Theorem \ref{theorem_disc_time_lower_bound}:} Let the
initial function be $\varphi =(\psi +\psi ^{\prime })/\sqrt{2}$ and let $%
\Delta :=d\left( \lambda ,\lambda ^{\prime }\right) .$ Then we compute:%
\begin{equation*}
d(T)=\frac{1}{2T}\left| 1-\cos \Delta T+\frac{\sin \Delta T\sin \Delta }{%
1-\cos \Delta }\right| Q\left( \psi ,\psi ^{\prime }\right) .
\end{equation*}

Consider interval $I_{r}=\Delta ^{-1}\pi (1/6+2r,5/6+2r),$ where $r$ is a
non-negative integer$.$ For every $T\in I_{r}$, $\sin \left( \Delta T\right)
\geq 1/2.$ Moreover, there is an integer $T_{r}\in I_{r}$ because $\left|
I_{r}\right| =\Delta ^{-1}\left( 2/3\right) \pi >1.$ Note that the distance
between $T_{r}$ and $T_{r+1}$ is less than $3\Delta ^{-1}\pi .$

For every integer $T_{r}$ we have 
\begin{equation*}
d(T_{r})\geq \frac{1}{4T_{r}}\left( \frac{\sin \Delta }{1-\cos \Delta }%
\right) Q\left( \psi ,\psi ^{\prime }\right) \geq \frac{1}{4T_{r}\Delta }%
Q\left( \psi ,\psi ^{\prime }\right) .
\end{equation*}%
(The second inequality holds because $\Delta \leq 2,$ and therefore $\sin
\Delta /\left( 1-\cos \Delta \right) \geq \Delta ^{-1}.$)

It follows that if $\varepsilon >0$ is sufficiently small then there exists
an integer $T_{r}$ between $Q\left( \psi ,\psi ^{\prime }\right)
/4\varepsilon \Delta -3\Delta ^{-1}\pi $ and $Q\left( \psi ,\psi ^{\prime
}\right) /4\varepsilon \Delta ,$ such that $d(T_{r})\geq \varepsilon .$ In
particular, if $\varepsilon \leq Q(\psi ,\psi ^{\prime })/80,$ then $%
T_{r}\geq Q\left( \psi ,\psi ^{\prime }\right) /4\varepsilon \Delta -3\Delta
^{-1}\pi \geq Q\left( \psi ,\psi ^{\prime }\right) /8\varepsilon \Delta ,$
and we conclude%
\begin{equation*}
t_{mix}\left( \varepsilon \right) \geq \frac{Q\left( \psi ,\psi ^{\prime
}\right) }{8\varepsilon \Delta }.
\end{equation*}%
QED.

\bigskip

\textbf{Proof of Theorem \ref{theorem_nonunitary_convergence}}: \ The proof
of (i) and (ii) is an application of results by Krein and Rutman from \cite%
{krein_rutman48}. In this paper a cone $K$ in a Banach space $X$ is fixed
and operator $A\in L\left( X\right) $ is called strongly positive if for
every non-zero $x\in K,$ there is an integer $n>0,$ such that $A^{n}x$ is in
the interior of $K$. Theorem 6.3 of this paper (on page 70 of the English
translation) shows that if $A$ is compact and strongly positive, then there
exists one and only one eigenvector of $A$ in the interior of $K$ and the
corresponding eigenvalue exceeds all others in absolute value. Moreover, the
proof of the theorem shows that this eigenvalue is simple. The claim of our
theorem follows if we apply the Krein-Rutman theorem to the cone of
positive-definite matrices. Indeed, by \ref{proposition_properties_channels}
there exists $\rho _{st}$ such that $\mathcal{T}\rho _{st}=\rho _{st}.$
Since $\mathcal{T}$ is strongly positive, hence $\rho _{st}$ is in the
interior of $K$ (i.e., strictly positive); by the Krein-Rutman theorem it is
the only eigenvector in the interior of $K,$ and its eigenvalue $1$ is simple%
$.$

For (iii), let $Z$ be the space of Hermitian matrices with zero trace, $%
Z=\left\{ \rho :\mathrm{tr}\left( \rho \right) =0\right\} .$ Then, $\mathcal{%
T}Z\subset Z$ and all the eigenvalues of $\mathcal{T}|_{Z}$ are less than $1$
in absolute value, because $1$ is a simple eigenvalue and $\rho _{st}\notin
Z.$ It follows that the spectral radius $\eta $ of $\mathcal{T}|_{Z}$ is
smaller than $1.$ Hence, $\lim_{n\rightarrow \infty }\left\| \left( \mathcal{%
T}|_{Z}\right) ^{n}\right\| ^{1/n}=\eta <1,$ which implies that $\mathcal{T}%
^{n}z\rightarrow 0$ for every $z\in Z.$ Since $\rho _{0}-\rho _{st}$ belongs
to $Z$ for every density matrix $\rho _{0},$ we conclude that $\mathcal{T}%
^{n}\rho _{0}\rightarrow \rho _{st}$ for every $\rho _{0}.$ QED.

\bibliographystyle{plain}
\bibliography{comtest}

\end{document}